\documentclass[amstex,12pt, amssymb]{article}

\usepackage{mathtext}
\usepackage[cp1251]{inputenc}
\usepackage[T2A]{fontenc}
\usepackage[dvips]{graphicx}
\usepackage{amsmath}
\usepackage{amssymb}
\usepackage{amsxtra}
\usepackage{latexsym}
\usepackage{ifthen}

\newcounter{lemma}[section]

\newcounter{corollary}[section]

\newcounter{remark}[section]

\newcounter{theorem}[section]

\newcounter{proposition}[section]

\numberwithin{equation}{section}

\def\XXint#1#2#3{{\setbox0=\hbox{$#1{#2#3}{\int}$}
     \vcenter{\hbox{$#2#3$}}\kern-.5\wd0}}

\def\cc{\setcounter{equation}{0}
\setcounter{figure}{0}\setcounter{table}{0}}

\begin{document}

\markboth{\centerline{V. Gutlyanskii, V. Ryazanov, A. Yefimushkin}}
{\centerline{On Hilbert, Riemann, Neumann and Poincare problems for
plane quasiregular mappings}}

\author{{V. Gutlyanskii, V. Ryazanov, A. Yefimushkin}}

\title{{\bf On Hilbert, Riemann, Neumann  and  Poincare\\ problems for plane quasiregular mappings
}}

\maketitle

\large \begin{abstract} Recall that the Hilbert (Riemann-Hilbert)
boundary value problem for the Beltrami equations was recently
solved for general settings in terms of nontangential limits and
principal asymptotic values. Here it is developed a new approach
making possible to obtain new results on tangential limits in
multiply connected domains. It is shown that the spaces of the found
solutions have the infinite dimension for prescribed families of
Jordan arcs terminating in almost every boundary point. We give also
applications of results obtained by us for the Beltrami equations to
the boundary value problems of Dirichlet, Riemann, Neumann and
Poincare for A-harmonic functions in the plane.

\medskip


\end{abstract}

\large \cc
\section{Introduction}

The boundary value problems of Dirichlet, Hilbert (Riemann-Hilbert),
Riemann, Neumann and Poincare are basic problems in the theory of
analytic functions which are closely interconnected, see e.g.
monographs \cite{G}, \cite{M} and \cite{V}, and also recent works
\cite{ER} and \cite{R1}-\cite{R*} for the history.

\medskip

Here we continue the development of the theory of these problems for
analytic functions as well as for the mappings with bounded
distortion by Reshetnyak, see e.g. \cite{Re}, either equivalently
for quasiregular mappings, see e.g. \cite{MRV}, or quasiconformal
functions, see e.g. \cite{LV}, i.e., regular solutions of the
Beltrami equations.

\medskip

Let $D$ be a domain in the complex plane $\mathbb C$ and $\mu:
D\to\mathbb C$ be a measurable function with $|\mu(z)|<1$ a.e.
Recall that a partial differential equation
\begin{equation}\label{1}
f_{\bar{z}}=\mu(z)\cdot f_z\
\end{equation}
where $f_{\bar z}={\bar\partial}f=(f_x+if_y)/2$, $f_{z}=\partial
f=(f_x-if_y)/2$, $z=x+iy$, $f_x$ and $f_y$ are partial derivatives
of the function $f$ in $x$ and $y$, respectively, is said to be a
{\bf Beltrami equation}. The Beltrami equation~\eqref{1} is said to
be {\bf nondegenerate} if $||\mu||_{\infty}<1$.

Note that analytic functions satisfy ~\eqref{1} with $\mu(z)\equiv
0$. Recall also that a {\bf regular solution} of the Beltrami
equation is a continuous, discrete and open mapping $f:{D}\to{\Bbb
C}$ of the Sobolev class $W_{\rm{loc}}^{1,1}$ satisfying~\eqref{1}
a.e. Note that, in the case of nondegenerate Beltrami equations, a
regular solution $f$ belongs to class $W_{\rm{loc}}^{1,p}$ for some
$p>2$ and, moreover, its Jacobian $J_f(z)\ne 0$ for almost all $z\in
D$, and it is called a {\bf quasiconformal function}, see e.g.
Chapter VI in \cite{LV}. Moreover, $f$ is called a {\bf
quasiconformal mapping} if in addition it is a homeomorphism.


It is clear that if we start to consider boundary value problems
with measurable boundary data, then the requests on the existence of
the limits at all points $\zeta\in\partial D$ and along all paths
terminating in $\zeta$ (as we{ll as the conception of the index)
lose any sense. Thus, the notion of solutions of such boundary value
problems should be widened. The nontangential limits from the
function theory of a complex variable were a suitable tool, see e.g.
\cite{ER} and \cite{R1}. In \cite{R*}, it was proposed an
alternative approach admitting tangential limits for analytic
functions. At present paper, we use this approach to the Beltrami
equations.


One of the relevant problems concerns to the measurement of sets on
boun\-da\-ries of domains. In this connection, note that the sets of
the length measure zero as well as of the harmonic measure zero are
invariant under conformal mappings, however, they are not invariant
under quasiconformal mappings as it follows from the famous
Ahlfors-Beurling example  of quasisymmetric mappings of the real
axis that are not absolutely continuous, see \cite{AB}. Hence we are
forced to apply here instead of them the so-called absolute harmonic
measure by Nevanlinna, in other words, logarithmic capacity whose
zero sets are invariant under quasiconformal mappings, see e.g.
\cite{N}.

\bigskip

\section{Definitions and preliminary remarks}

By the well--known Priwalow uniqueness theorem analytic functions in
the unit disk $\mathbb D = \{ z\in\mathbb{C}: |z|<1\}$ coincide if
they have the equal boundary values along all nontangential paths to
a set $E$ of points in $\partial\mathbb D$ of a positive measure,
see e.g. Theorem IV.2.5 in \cite{P}. The theorem is valid also for
analytic functions in Jordan domains with rectifiable boundaries,
see e.g. Section IV.2.6 in \cite{P}.

\medskip

However, examples of Lusin and Priwalow show that there exist
nontrivial analytic functions in $\mathbb D$ whose radial boundary
values are equal to zero on sets $E\subseteq\partial\mathbb D$ of a
positive measure, see e.g. Section IV.5 in \cite{P}. Simultaneously,
by Theorem IV.6.2 in \cite{P} of Lusin and Priwalow the uniqueness
result is valid if $E$ is of the second category. Theorem 1 in
\cite{BS} demonstrates that the latter condition is necessary.

\medskip

Recall Baire's terminology for categories of sets and functions.
Namely, given a topological space $X$, a set $E\subseteq X$ is of
{\bf first category} if it can be written as a countable union of
nowhere dense sets, and is of {\bf second category} if $E$ is not of
first category. Also, given topological spaces $X$ and $X_*$,
$f:X\to X_*$ is said to be a {\bf function of Baire class 1} if
$f^{-1}(U)$ for every open set $U$ in $X_*$ is an $F_{\sigma}$ set
in $X$ where an $F_{\sigma}$ set is the union of a sequence of
closed sets.

\medskip

Given a Jordan curve $C$ in $\mathbb C$, we say that a family of
Jordan arcs $\{J_{\zeta}\}_{\zeta\in C}$ is of {\bf class}
${\cal{BS}}$ {\bf (of the Bagemihl--Seidel class)}, cf. \cite{BS},
740--741, if all $J_{\zeta}$ lie in a ring $\frak{R}$ generated by
$C$ and a Jordan curve $C_*$ in $\mathbb C$, $C_*\cap C={\O}$,
$J_{\zeta}$ is joining $C_*$ and $\zeta\in C$, every $z\in \frak{R}$
belongs to a single arc $J_{\zeta}$, and for a sequence of mutually
disjoint Jordan curves $C_n$ in $\frak{R}$ such that $C_n\to C$ as
$n\to\infty$, $J_{\zeta}\cap C_n$ consists of a single point for
each $\zeta\in C$ and $n=1,2, \ldots$.

\medskip

In particular, a family of Jordan arcs $\{J_{\zeta}\}_{\zeta\in C}$
is of class ${\cal{BS}}$ if $J_{\zeta}$ are in\-du\-ced by an
isotopy of $C$. For instance, every curvilinear ring $\frak{R}$ one
of whose boundary component is $C$ can be mapped with a conformal
mapping $g$ onto a circular ring $R$ and the inverse mapping
$g^{-1}:R\to\frak{R}$ maps radial lines in $R$ onto suitable Jordan
arcs $J_{\zeta}$ and centered circles in $R$ onto Jordan curves
giving the corresponding isotopy of $C$ to other boundary component
of $\frak{R}$.

\bigskip

Now, if $\Omega\subset\mathbb C$ is an open set bounded by a finite
collection of mutually disjoint Jordan curves, then we say that a
family of Jordan arcs $\{J_{\zeta}\}_{\zeta\in \partial\Omega}$ is
of class ${\cal{BS}}$ if its restriction to each component of
$\partial\Omega$ is so.

\medskip

Theorem 1 in \cite{BS} can be written in the following way, see
footnote 9 in \cite{BS}.

\medskip

{\bf Proposition 1.} {\it  Let $D$ be a bounded domain in $\mathbb
C$ whose boundary consists of a finite number of mutually disjoint
Jordan curves and let $\{ \gamma_{\zeta}\}_{\zeta\in\partial D}$ be
a family of Jordan arcs of class ${\cal{BS}}$ in ${D}$.


Suppose $M$ is an $F_{\sigma}$ set of first category on
$\partial{D}$ and $\Phi(\zeta)$ is a complex-valued function of
Baire class 1 on $M$. Then there is a nonconstant single-valued
analytic function $f: D\to\mathbb C$ such that,  for all $\zeta\in
M$, along $\gamma_{\zeta}$
\begin{equation}\label{eqBS} \lim\limits_{z\to\zeta}\ f(z)\ =\
\Phi(\zeta)\ .
\end{equation} }

On this basis, in the case of domains $D$ whose boundaries consist
of rectifiable Jordan curves, it was formulated Theorem 2 in
\cite{BS} on the existence of analytic functions $f: D\to\mathbb C$
such that (\ref{eqBS}) holds a.e. on $\partial{D}$ with respect to
the natural parameter for each prescribed measurable function
$\Phi:\partial{D}\to\mathbb C$.

\bigskip

\section{On the Dirichlet problem for analytic functions}

\bigskip

We need for our goals the following statement on analytic functions
with prescribed boundary values that is similar to Theorem 2 in
\cite{BS} but formulated in terms of logarithmic capacity instead of
the natural parameter, see definitions, notations and comments in
\cite{ER}.

\bigskip

{\bf Theorem 1.} {\it  Let $D$ be a bounded domain in $\mathbb C$
whose boundary consists of a finite number of mutually disjoint
Jordan curves and let a function $\Phi:\partial D\to\mathbb C$ be
measurable with respect to the logarithmic capacity.

\medskip
Suppose that $\{ \gamma_{\zeta}\}_{\zeta\in\partial D}$ is a family
of Jordan arcs of class ${\cal{BS}}$ in ${D}$. Then there is a
nonconstant single-valued analytic function $f: D\to\mathbb C$ such
that (\ref{eqBS}) holds
along $\gamma_{\zeta}$ for a.e. $\zeta\in\partial D$ with respect to
the logarithmic capacity. }

\bigskip

\begin{proof}
Note first of all that ${\cal{C}}:=C(\partial D)<\infty$ because
$\partial D$ is bounded and Borel, even compact, and show that there
is a sigma-compact set $S$ in $\partial D$ of first category such
that $C(S)={\cal{C}}$. More precisely, $S$ will be the union of a
sequence of sets $S_m$ in $\partial D$ of the Cantor type that are
nowhere dense in $\partial D$.

Namely, $S_m$ is constructed in the following way. First we remove
an open arc $A_1$ in $\partial D$ of the logarithmic capacity
$2^{-m}{\cal{C}}$ and one more open arc $A_2$ in $\partial
D\setminus A_1$ of the logarithmic capacity $2^{-(m+1)}{\cal{C}}$
such that $\partial D\setminus (A_1\cup A_2)$ consists of 2 segments
of $\partial D$ with the equal logarithmic capacity. Then we remove
a union $A_3$ of 2 open arcs in each these segments of the total
logarithmic capacity $2^{-(m+2)}{\cal{C}}$ such that new 4 segments
in $\partial D\setminus (A_1\cup A_2\cup A_3)$ have the equal
logarithmic capacity. Repeating by induction this construction, we
obtain the compact sets $S_m=\partial
D\setminus\bigcup\limits_{m=1}\limits^{\infty } A_m$ with the
logarithmic capacity $(1-2^{-(m-1)})\cdot{\cal{C}}\to{\cal{C}}$ as
$m\to\infty$.

Note also that the logarithmic capacity is Borel's regular measure
as well as Radon's measure in the sense of points 2.2.3 and 2.2.5 in
\cite{Fe}, correspondingly, see Section 2 in \cite{ER}. Hence the
Lusin theorem holds for the lo\-ga\-rith\-mic capacity on $\mathbb
C$, see Theorem 2.3.5 in \cite{Fe}.

By the Lusin theorem one can find a sequence of compacta $K_n$ in
$S$ with $C(S\setminus K_n)<2^{-n}$  such that $\Phi|_{K_n}$ is
continuous for each $n=1,2,\ldots$, i.e., for every open set
$U\subseteq\mathbb C$, $W_n:=\Phi|_{K_n}^{-1}(U)=V_n\cap S$ for some
open set $V_n$ in $\mathbb C$, and $W_n$ is sigma-compact because
$V_n$ and $S$ are sigma-compact. Consequently, $W:=\Phi|_{K}^{-1}(U)
=\bigcup\limits_{n=1}\limits^{\infty } W_n$ is also sigma-compact
where $K = \bigcup\limits_{n=1}\limits^{\infty } K_n$. Hence the
restriction of $\Phi$ on the set $K$ is a function of Baire's class
1. Finally, note that by the construction $C(Z)=0$ where $Z=\partial
D\setminus K=(\partial D\setminus S)\cup (S\setminus K)$ and
$$
K = \bigcup\limits_{m,\, n=1}\limits^{\infty } (S_m\cap K_n)
$$
where each set $E_{m,n}:=S_m\cap K_n$, $m,n=1,2,\ldots$, is nowhere
dense in $\partial D$.

\medskip

Thus, the conclusion of Theorem 1 follows from Proposition 1.
\end{proof}\ $\Box$

\bigskip

{\bf Corollary 1.} {\it  Let $D$ be a bounded domain in $\mathbb C$
whose boundary consists of a finite number of mutually disjoint
Jordan curves and let a function $\varphi:\partial D\to\mathbb R$ be
measurable with respect to the logarithmic capacity.

\medskip
Suppose that $\{ \gamma_{\zeta}\}_{\zeta\in\partial D}$ is a family
of Jordan arcs of class ${\cal{BS}}$ in ${D}$. Then there is a
harmonic function $u: D\to\mathbb R$ such that
\begin{equation}\label{eqHARMONIC} \lim\limits_{z\to\zeta}\ u(z)\ =\ \varphi(\zeta)
\end{equation}
along $\gamma_{\zeta}$ for a.e. $\zeta\in\partial D$ with respect to
the logarithmic capacity. }

\bigskip

{\bf Corollary 2.} {\it  Let $D$ be a Jordan domain in $\mathbb C$
and a function $\Phi:\partial D\to\mathbb C$ be measurable with
respect to the logarithmic capacity.

\medskip
Suppose that $\{ \gamma_{\zeta}\}_{\zeta\in\partial D}$ is a family
of Jordan arcs of class ${\cal{BS}}$ in ${D}$. Then there is a
harmonic function $u: D\to\mathbb R$ such that
\begin{equation}\label{eqGRAD} \lim\limits_{z\to\zeta}\ \nabla u(z)\ =\ \Phi(\zeta)
\end{equation}
along $\gamma_{\zeta}$ for a.e. $\zeta\in\partial D$ with respect to
the logarithmic capacity. }

\medskip

Here we use the complex writing for the gradient ${\nabla u} := u_x
+ i\cdot u_y$.

\medskip

\begin{proof}
Indeed, by Theorem 1 there is a single-valued analytic function $f:
D\to\mathbb R$ such that \begin{equation}\label{eqGRADCON}
\lim\limits_{z\to\zeta}\ f(z)\ =\ \overline{\Phi(\zeta)}
\end{equation} along $\gamma_{\zeta}$ for a.e.
$\zeta\in\partial D$ with respect to the logarithmic capacity. Then
any indefinite integral $F$ of $f$ is also a single-valued analytic
function in the simply connected domain $D$ and the harmonic
functions $u=\mathrm {Re}\ F$ and $v=\mathrm {Im}\ F$  satisfy the
Cauchy-Riemann system $v_x=-u_y$ and $v_y=u_x$. Hence
$$
f\ =\ F^{\prime}\ =\ F_x\ =\ u_x\ +\ i\cdot v_x\ =\ \ u_x\ -\ i\cdot
u_y\ =\ \overline{\nabla u}\ .
$$
Thus, (\ref{eqGRAD}) follows from (\ref{eqGRADCON}) and,
consequently, $u$ is the desired function.
\end{proof}\ $\Box$

\bigskip

{\bf Corollary 3.} {\it  Let a function $\Phi:\partial \mathbb
D\to\mathbb C$ be measurable with respect to the logarithmic
capacity. Then there is a harmonic function $u:\mathbb D\to\mathbb
R$ such that (\ref{eqGRAD}) holds along radial lines for a.e.
$\zeta\in\partial\mathbb D$ with respect to the logarithmic
capacity. }

\bigskip

\section{On the Dirichlet problem for quasiregular mappings}

We show that quasiconformal functions are able to take arbitrary
measurable values along prescribed families of Jordan arcs
terminating in the boundary.


{\bf Theorem 2.} {\it  Let $D$ be a bounded domain in $\mathbb C$
whose boundary consists of a finite number of mutually disjoint
Jordan curves, $\mu:D\to\mathbb{C}$ be a measurable (by Lebesgue)
function with $||\mu||_{\infty}<1$, and let a function
$\Phi:\partial D\to\mathbb C$ be measurable with respect to the
logarithmic capacity.

Suppose that $\{ \gamma_{\zeta}\}_{\zeta\in\partial D}$ is a family
of Jordan arcs of class ${\cal{BS}}$ in ${D}$. Then the Beltrami
equation \eqref{1} has a regular solution $f: D\to\mathbb C$ such
that (\ref{eqBS}) holds
along $\gamma_{\zeta}$ for a.e. $\zeta\in\partial D$ with respect to
the logarithmic capacity. }


\begin{proof}
Extending $\mu$ by zero everywhere outside of $D$, we obtain the
existence of a quasiconformal mapping
$h:{\mathbb{C}}\to{\mathbb{C}}$ a.e. satisfying the Beltrami
equation \eqref{1} with the given $\mu$, see e.g. Theorem V.B.3 in
\cite{Alf}. Setting $D_*=h(D)$ and
$\Gamma_{\xi}=h(\gamma_{h^{-1}(\xi)})$, $\xi\in\partial D_*$, we see
that $\{ \Gamma_{\xi}\}_{\xi\in\partial D_*}$ is a family of Jordan
arcs of class ${\cal{BS}}$ in $D_*$.

The logarithmic capacity of a set coincides with its transfinite
diameter, see e.g.  \cite{F} and the point 110 in \cite{N}. Hence
the mappings $h$ and $h^{-1}$ transform sets of logarithmic capacity
zero on $\partial D$ into sets of logarithmic capacity zero on
$\partial D_*$ and vice versa because these quasiconformal mappings
are continuous by H\"older on $\partial D$ and $\partial D_*$
correspondingly, see e.g. Theorem II.4.3 in \cite{LV}.

Further, the function $\varphi=\Phi\circ h^{-1}$ is measurable with
respect to logarithmic capacity. Indeed, under this mapping
measurable sets with respect to lo\-ga\-rith\-mic capacity are
transformed into measurable sets with respect to logarithmic
capacity because such a set can be represented as the union of a
sigma-com\-pac\-tum and a set of logarithmic capacity zero and
compacta under continuous mappings are transformed into compacta and
compacta are measurable sets with respect to logarithmic capacity.

By Theorem 1 there is an analytic function $F: D_*\to\mathbb C$ such
that
\begin{equation}\label{29}
\lim_{w\to\xi}\ F(w)\ =\ \varphi(\xi)
\end{equation}
holds along $\Gamma_{\xi}$ for a.e. $\xi\in\partial D_*$ with
respect to logarithmic capacity. Thus, by the arguments given above
$f:=F\circ h$ is the desired solution of \eqref{1}.
\end{proof}\ $\Box$

\bigskip

\section{On Riemann-Hilbert problem for analytic functions}

\medskip

Recall that the classical setting of the {\bf Hilbert
(Riemann-Hilbert) boun\-da\-ry value problem} was on finding
analytic functions $f$ in a domain $D\subset\mathbb C$ bounded by a
rectifiable Jordan curve with the boundary condition
\begin{equation}\label{eqHILBERT}
\lim\limits_{z\to\zeta}\ \mathrm {Re}\
\{\overline{\lambda(\zeta)}\cdot f(z)\}\ =\ \varphi(\zeta)
\quad\quad\quad\ \ \ \forall \ \zeta\in \partial D
\end{equation}
where functions $\lambda$ and $\varphi$ were continuously
differentiable with respect to the natural parameter $s$ on
$\partial D$ and, moreover, $|\lambda|\ne 0$ everywhere on $\partial
D$. Hence without loss of generality one can assume that
$|\lambda|\equiv 1$ on $\partial D$.

\medskip

{\bf Theorem 3.} {\it\, Let $D$ be a bounded domain in $\mathbb C$
whose boundary consists of a finite number of mutually disjoint
Jordan curves, and let $\lambda:\partial D\to\mathbb C$, $|\lambda
(\zeta)|\equiv 1$, $\varphi:\partial D\to\mathbb R$ and
$\psi:\partial D\to \mathbb R$ be measurable functions with respect
to the logarithmic capacity.

\medskip

Suppose that $\{ \gamma_{\zeta}\}_{\zeta\in\partial D}$ is a family
of Jordan arcs of class ${\cal{BS}}$ in ${D}$. Then there is a
nonconstant single-valued analytic function $f: D\to\mathbb C$ such
that
\begin{equation}\label{eqARE} \lim\limits_{z\to\zeta}\ \mathrm {Re}\
\{\overline{\lambda(\zeta)}\cdot f(z)\}\ =\ \varphi(\zeta)\ ,
\end{equation}
\begin{equation}\label{eqIM}
\lim\limits_{z\to\zeta}\ \mathrm {Im}\
\{\overline{\lambda(\zeta)}\cdot f(z)\}\ =\ \psi(\zeta)\
\end{equation}
along $\gamma_{\zeta}$ for a.e. $\zeta\in\partial D$ with respect to
the logarithmic capacity.}

\medskip

{\bf Remak 1.} Thus, the space of all solutions $f$ of the Hilbert
problem (\ref{eqARE}) in the given sense has the infinite dimension
for any prescribed $\varphi$, $\lambda$ and $\{
\gamma_{\zeta}\}_{\zeta\in D}$ because the space of all functions
$\psi:\partial D\to \mathbb R$ which are  measurable with respect to
the logarithmic capacity has the infinite dimension.

The latter is valid even for its subspace of continuous functions
$\psi:\partial D\to \mathbb R$. Indeed, by the Riemann theorem every
Jordan domain $G$ can be mapped with a conformal mapping $g$ onto
the unit disk $\mathbb D$ and by the Caratheodory theorem $g$ can be
extended to a homeomorphism of $\overline G$ onto $\overline{\mathbb
D}$. By the Fourier theory, the space of all continuous functions
$\tilde\psi:\partial\mathbb D\to\mathbb R$, equivalently, the space
of all continuous $2\pi$-periodic functions $\psi_*:\mathbb
R\to\mathbb R$, has the infinite dimension.

\medskip

\begin{proof} Indeed, set $\Psi(\zeta)=\varphi(\zeta)+i\cdot\psi(\zeta)$
and $\Phi(\zeta)=\lambda(\zeta)\cdot\Psi(\zeta)$ for all
$\zeta\in\partial D$. Then by Theorem 1 there is a single-valued
analytic function $f$ such that
\begin{equation}\label{eqABS} \lim\limits_{z\to\zeta}\ f(z)\ =\ \Phi(\zeta)
\end{equation}
along $\gamma_{\zeta}$ for a.e. $\zeta\in\partial D$ with respect to
the logarithmic capacity. Then also \begin{equation}\label{eqLIM}
\lim\limits_{z\to\zeta}\ \overline{\lambda(\zeta)}\cdot f(z)\ =\
\Psi(\zeta)
\end{equation}
along $\gamma_{\zeta}$ for a.e. $\zeta\in\partial D$ with respect to
the logarithmic capacity.
\end{proof}\ $\Box$

\bigskip

\section{On Riemann-Hilbert problem for Beltrami equations}

{\bf Theorem 4.} {\it\, Let $D$ be a bounded domain in $\mathbb C$
whose boundary consists of a finite number of mutually disjoint
Jordan curves, $\mu:D\to\mathbb{C}$ be a measurable (by Lebesgue)
function with $||\mu||_{\infty}<1$, and let $\lambda:\partial
D\to\mathbb C$, $|\lambda (\zeta)|\equiv 1$, $\varphi:\partial
D\to\mathbb R$ and $\psi:\partial D\to \mathbb R$ be functions that
are  measurable with respect to the logarithmic capacity.

\medskip

Suppose that $\{ \gamma_{\zeta}\}_{\zeta\in\partial D}$ is a family
of Jordan arcs of class ${\cal{BS}}$ in ${D}$. Then the Beltrami
equation \eqref{1} has a regular solution $f: D\to\mathbb C$ such
that
\begin{equation}\label{eqBRE} \lim\limits_{z\to\zeta}\ \mathrm {Re}\
\{\overline{\lambda(\zeta)}\cdot f(z)\}\ =\ \varphi(\zeta)\ ,
\end{equation}
\begin{equation}\label{eqBIM}
\lim\limits_{z\to\zeta}\ \mathrm {Im}\
\{\overline{\lambda(\zeta)}\cdot f(z)\}\ =\ \psi(\zeta)\
\end{equation}
along $\gamma_{\zeta}$ for a.e. $\zeta\in\partial D$ with respect to
the logarithmic capacity.}

\medskip

\begin{proof} Indeed, set $\Psi(\zeta)=\varphi(\zeta)+i\cdot\psi(\zeta)$
and $\Phi(\zeta)=\lambda(\zeta)\cdot\Psi(\zeta)$ for all
$\zeta\in\partial D$. Then by Theorem 2 the Beltrami equation
\eqref{1} has a regular solution $f: D\to\mathbb C$ such that
\begin{equation}\label{eqBBS} \lim\limits_{z\to\zeta}\ f(z)\ =\ \Phi(\zeta)
\end{equation}
along $\gamma_{\zeta}$ for a.e. $\zeta\in\partial D$ with respect to
the logarithmic capacity. Then also \begin{equation}\label{eqBLIM}
\lim\limits_{z\to\zeta}\ \overline{\lambda(\zeta)}\cdot f(z)\ =\
\Psi(\zeta)
\end{equation}
along $\gamma_{\zeta}$ for a.e. $\zeta\in\partial D$ with respect to
the logarithmic capacity.
\end{proof}\ $\Box$

\medskip

{\bf Remak 2.} Thus, the space of all solutions $f$ of the Hilbert
problem (\ref{eqBRE}) for  the Beltrami equation \eqref{1} in the
given sense has the infinite dimension for any prescribed $\varphi$,
$\lambda$ and $\{ \gamma_{\zeta}\}_{\zeta\in D}$ because the space
of all functions $\psi:\partial D\to \mathbb R$ which are measurable
with respect to the logarithmic capacity has the infinite dimension,
see Remark 1.

\bigskip

\section{Neumann and Poincare problems for analytic functions}

As known the Neumann problem has no classical solutions generally
spea\-king even for harmonic functions and smooth boundary data, see
e.g. \cite{Mi}. Here we study it and more general boundary value
problems for analytic functions with measurable boundary data in a
generalized form.

\medskip

First of all recall a more general {\bf problem on directional
derivatives} for the harmonic functions in the unit disk $\mathbb D
= \{ z\in\mathbb{C}: |z|<1\}$, $z=x+iy$. The classic setting of the
latter problem is to find a function $u:\mathbb D\to\mathbb R$ that
is twice continuously differentiable, admits a continuous extension
to the boundary of $\mathbb D$ together with its first partial
derivatives, satisfies the Laplace equation
\begin{equation}\label{eqLAPLACE}
\Delta u\ :=\ \frac{\partial^2 u}{\partial x^2}\ +\ \frac{\partial^2
u}{\partial y^2}\ =\ 0 \quad\quad\quad\forall\ z\in\mathbb D
\end{equation} and the boundary condition with a prescribed
continuous date $\varphi : \partial\mathbb D\to\mathbb R$:
\begin{equation}\label{eqDIRECT}
\frac{\partial u}{\partial \nu}\ =\ \varphi(\zeta) \quad\quad\quad
\forall\ \zeta\in\partial\mathbb D
\end{equation}
where $\frac{\partial u}{\partial \nu}$ denotes the derivative of
$u$ at $\zeta$ in a direction $\nu = \nu(\zeta)$, $|\nu(\zeta)|=1$:
\begin{equation}\label{eqDERIVATIVE}
\frac{\partial u}{\partial \nu}\ :=\ \lim_{t\to 0}\
\frac{u(\zeta+t\cdot\nu)-u(\zeta)}{t}\ .
\end{equation}

{\bf The Neumann problem} is a special case of the above problem on
directional derivatives with  the boundary condition
\begin{equation}\label{eqNEUMANN}
\frac{\partial u}{\partial n}\ =\ \varphi(\zeta) \quad\quad\quad
\forall\ \zeta\in\partial\mathbb D
\end{equation}
where $n$ denotes the unit interior normal to $\partial\mathbb D$ at
the point $\zeta$.


In turn, the above problem on directional derivatives is a special
case of {\bf the Poincare problem} with  the boundary condition
\begin{equation}\label{eqPOINCARE}
a\cdot u\ +\ b\cdot\frac{\partial u}{\partial \nu}\ =\
\varphi(\zeta) \quad\quad\quad \forall\ \zeta\in\partial\mathbb D
\end{equation}
where $a=a(\zeta)$ and $b=b(\zeta)$ are real-valued functions given
on $\partial\mathbb D$.


Recall also that twice continuously differentiable solutions of the
Laplace equation are called {\bf harmonic functions}. As known, such
functions are infinitely differentiable as well as they are real and
imaginary parts of analytic functions. The above boundary value
problems for analytic functions are formulated in a similar way. Let
us start from the problem on directional derivatives.

\medskip

{\bf Theorem 5.} {\it\, Let $D$ be a Jordan domain in $\mathbb C$
with a rectifiable boundary, $\nu:\partial D\to\mathbb C$, $|\nu
(\zeta)|\equiv 1$, and $\Phi:\partial D\to\mathbb C$ be measurable
functions with respect to the natural parameter.

Suppose that $\{ \gamma_{\zeta}\}_{\zeta\in\partial D}$ is a family
of Jordan arcs of class ${\cal{BS}}$ in ${D}$. Then there is a
nonconstant single-valued analytic function $f: D\to\mathbb C$ such
that
\begin{equation}\label{eqDLIMIT-R} \lim\limits_{z\to\zeta}\ \frac{\partial f}{\partial \nu}\ (z)\ =\
\Phi(\zeta)\end{equation} along $\gamma_{\zeta}$ for a.e.
$\zeta\in\partial D$ with respect to the natural parameter.}

\bigskip

\begin{proof} Indeed, by Theorem 5 in \cite{BS} there is a harmonic function $u: D\to\mathbb R$ such that
\begin{equation}\label{eqGRADIENT} \lim\limits_{z\to\zeta} \nabla u(z)\ =\ \Psi(\zeta)\ :=\ {\nu(\zeta)}\cdot \overline{\Phi(\zeta)}\end{equation}
along $\gamma_{\zeta}$ for a.e. $\zeta\in\partial D$ with respect to
the natural parameter where we use the complex writing for the
gradient $\nabla u = u_x+i\cdot u_y\ $.

As known, there is a conjugate harmonic function $v$ such that
$f:=u+i\cdot v$ is a single-valued analytic function in the simply
connected domain $D$, see e.g. point I.A in \cite{Ku}. The functions
$u$ and $v$ satisfy the Cauchy-Riemann system $v_x=-u_y$ and
$v_y=u_x$, i.e., $\nabla v=i\cdot\nabla u$. Note also that the
derivatives of $u$ and $v$ in the direction $\nu$ are projections of
gradients of $u$ and $v$ onto $\nu$. Hence
$$
\frac{\partial f}{\partial \nu}\ =\ \frac{\partial u}{\partial \nu}
+ i\cdot\frac{\partial v}{\partial \nu}\ =\ \mathrm {Re}\
\nu\cdot\overline{\nabla u}\ +\ i\cdot\mathrm {Re}\
\nu\cdot\overline{\nabla v}\ =\ {\nu\cdot\overline{\nabla u}}\ .
$$
Thus, (\ref{eqDLIMIT-R}) follows from (\ref{eqGRADIENT}) and,
consequently, $f$ is the desired function.
\end{proof}\ $\Box$

\bigskip

{\bf Remark 3.} We are able to say more in the case $\mathrm {Re}\
n\cdot\overline{\nu}>0$ where $n=n(\zeta)$ is the unit interior
normal with a tangent to $\partial D$ at a point $\zeta\in\partial
D$. In view of (\ref{eqDLIMIT-R}), since the limit $\Phi(\zeta)$ is
finite, there is a finite limit $f(\zeta)$ of $f(z)$ as $z\to\zeta$
in $D$ along the straight line passing through the point $\zeta$ and
being parallel to the vector $\nu(\zeta)$ because along this line,
for $z$ and $z_0$ that are close enough to $\zeta$,
\begin{equation}\label{eqDIFFERENCE-R} f(z)\ =\ f(z_0)\ -\ \int\limits_{0}\limits^{1}\
\frac{\partial f}{\partial \nu}\ (z_0+\tau (z-z_0))\ d\tau\
.\end{equation} Thus, at each point with the condition
(\ref{eqDLIMIT-R}), there is the directional derivative
\begin{equation}\label{eqPOSITIVE-R}
\frac{\partial f}{\partial \nu}\ (\zeta)\ :=\ \lim_{t\to 0}\
\frac{f(\zeta+t\cdot\nu)-f(\zeta)}{t}\ =\ \Phi(\zeta)\ .
\end{equation}

\medskip

In particular, we obtain on the basis of Theorem 5 and Remark 3 the
following statement on the Neumann problem.

\medskip

{\bf Corollary 4.} {\it\, For each measurable function
$\Phi:\partial\mathbb D\to\mathbb C$, one can find a nonconstant
single-valued analytic function $f:\mathbb D\to\mathbb C$ such that,
for a.e. point $\zeta\in\partial\mathbb D$, there exist:

\bigskip

1) the finite radial limit
\begin{equation}\label{eqLIMIT1}
f(\zeta)\ :=\ \lim\limits_{r\to 1}\ f(r\zeta)\end{equation}

2) the normal derivative
\begin{equation}\label{eqNORMAL}
\frac{\partial f}{\partial n}\ (\zeta)\ :=\ \lim_{t\to 0}\
\frac{f(\zeta+t\cdot n)-f(\zeta)}{t}\ =\ \Phi(\zeta)
\end{equation}

3) the radial limit
\begin{equation}\label{eqLIMIT2} \lim\limits_{r\to 1}\ \frac{\partial f}{\partial n}\ (r\zeta)\ =\
\frac{\partial f}{\partial n}\ (\zeta)\end{equation} where
$n=n(\zeta)$ denotes the unit interior normal to $\partial\mathbb D$
at the point $\zeta$.}

\bigskip

Finally, arguing in a perfectly similar way to the proof of Theorem
5 but basing on Corollary 2 above instead of Theorem 5 in \cite{BS},
we obtain the following theorem that solves the problem on
directional derivatives for analytic functions in arbitrary Jordan
domains.

\medskip

{\bf Theorem 6.} {\it\, Let $D$ be a Jordan domain in $\mathbb C$,
$\nu:\partial D\to\mathbb C$, $|\nu (\zeta)|\equiv 1$, and
$\Phi:\partial D\to\mathbb C$ be measurable functions with respect
to the logarithmic capacity.

Suppose that $\{ \gamma_{\zeta}\}_{\zeta\in\partial D}$ is a family
of Jordan arcs of class ${\cal{BS}}$ in ${D}$. Then there is a
nonconstant single-valued analytic function $f: D\to\mathbb C$ such
that
\begin{equation}\label{eqA-NP-D-LIMIT-R} \lim\limits_{z\to\zeta}\ \frac{\partial f}{\partial \nu}\ (z)\ =\
\Phi(\zeta)\end{equation} along $\gamma_{\zeta}$ for a.e.
$\zeta\in\partial D$ with respect to the logarithmic capacity.}

\medskip

{\bf Remark 4.} As it follows from the given scheme of the proof,
Theorems 5 and 6 are valid for multiply connected domains bounded by
a finite collection of mutually disjoint Jordan curves, however, the
corresponding analytic function $f$ can be multivalent.

\bigskip

\section{Neumann, Poincare problems for Beltrami equations}

{\bf Theorem 7.} {\it\, Let $D$ be a Jordan domain in $\mathbb C$,
$\mu:D\to\mathbb{C}$ be a function of the  H\"older class
$C^{\alpha}$ with $\alpha\in(0,1)$ and $|\mu(z)|\leq k<1$, $z\in D$,
and let $\nu:\partial D\to\mathbb C$, $|\nu (\zeta)|\equiv 1$, and
$\Phi:\partial D\to\mathbb C$ be measurable with respect to the
logarithmic capacity.

\medskip

Suppose that $\{ \gamma_{\zeta}\}_{\zeta\in\partial D}$ is a family
of Jordan arcs of class ${\cal{BS}}$ in ${D}$. Then the Beltrami
equation \eqref{1} has a regular solution $f: D\to\mathbb C$ of the
class $C^{1+\alpha}$ such that
\begin{equation}\label{eqNP-LIMIT-R} \lim\limits_{z\to\zeta}\ \frac{\partial f}{\partial \nu}\ (z)\ =\
\Phi(\zeta)\end{equation} along $\gamma_{\zeta}$ for a.e.
$\zeta\in\partial D$ with respect to the logarithmic capacity.}




{\bf Lemma 1.} {\it Let $S$ be a subspace of a metric space $(X,d)$
and let $h: S\to\mathbb R$ be a H\"older continuous function,
namely, for some $C>0$ and $\alpha\in(0,1)$,
\begin{equation}\label{eqNP1}
|h(x_1)-h(x_2)|\ \leq\ C\cdot d^{\alpha}(x_1,x_2)\ \ \ \ \ \
\forall\  x_1,\ x_2\in S\ .
\end{equation}
Then there is an extension $H:X\to\mathbb R$ of $h$ to $X$ such that
\begin{equation}\label{eqNP2}
|H(x_1)-H(x_2)|\ \leq\ C\cdot d^{\alpha}(x_1,x_2)\ \ \ \ \ \
\forall\ x_1,\ x_2\in X\ .
\end{equation} The same is valid for the functions $h:
S\to\mathbb C$ as well as $h: S\to\mathbb R^n$.}

\begin{proof}
Let us prove that $d_*(x_1,x_2):=d^{\alpha}(x_1,x_2)$,
$(x_1,x_2)\in(X\times X)$, is a distance on $X$. Indeed, by the
triangle inequality for $d$ we have that
$$
d_*(x_1,x_2)\ \leq\ (d(x_1,x_3)\ +\ d(x_3,x_2))^{\alpha}\ \ \ \ \ \
\forall\ x_1,\ x_2,\ x_3\in X\ .
$$
Hence it remains to show that
$$
\varphi(t,\tau)\ :=\ t^{\alpha}\ +\ \tau^{\alpha}\ -\ (t\ +\
\tau)^{\alpha}\ \geq\ 0\ \ \ \ \ \ \ \forall\ t,\tau\in\mathbb R^+,\
\alpha\in(0,1)\ .
$$
It sufficient for this to note that $\varphi(t,0)\equiv 0$ and
$$
\frac{\partial \varphi}{\partial \tau}\ (t,\tau)\ =\
\alpha\cdot\left[ \frac{1}{\tau^{1-\alpha}}\ -\ \frac{1}{(t\ +\
\tau)^{1-\alpha}}\right]\ \geq\ 0\ .
$$

Thus, the conclusion of Lemma 1 follows from point 2.10.44 in
\cite{Fe} on the Lipschitz functions.
\end{proof}\ $\Box$

\bigskip

{\bf Lemma 2.} {\it Let $f:D\to\mathbb C$ be a smooth quasiconformal
mapping satisfying the Beltrami equation \eqref{1} with a continuous
complex coefficient $\mu$ in a domain $D\subseteq\mathbb C$. Then
\begin{equation}\label{eqNP3}
C(Z_f)\ =\ 0\ =\ C(f(Z_f))\ ,
\end{equation}
where $C(E)$ is the logarithmic capacity of $E\subseteq\mathbb C$,
$J_f$ is the Jacobian of $f$, and
\begin{equation}\label{eqNP4}
Z_f\ =\ \{\ z\in\mathbb C: J_f(z)\ =\ 0\ \}\ .
\end{equation}}

\begin{proof}
First of all, note that $f$ satisfies \eqref{1} for all $z\in D$
because the functions $\mu(z)$, $f_z$ and $f_{\bar z}$ are
continuous, and by \eqref{1} the maximal distortion
$M_f(z)=|f_z|+|f_{\bar z}|$ under $f$ at $z$ is equal to $0$ at all
points of $Z_f$. Note also that, in view of countable subadditivity
of logarithmic capacity, with no loss of generality we may assume
that $D$ is bounded and then $C(D)<\infty$.

 Let us take an arbitrary $\varepsilon>0$. In view of
continuity of $M_f(z)$, for every $z_*\in Z_f$, there is a disk
$B_{z_*}$ centered at $z_*$ with a small enough radius such that
${B_{z_*}}\subset D$ and such that $\sup\limits_{z\in {B_{z_*}}}\,
M_f(z)\leq\varepsilon$. The union of all these disks gives an open
neighborhood $U$ of $Z_f$ such that $\sup\limits_{z\in U}\,
M_f(z)\leq\varepsilon$. The logarithmic capacity of a set can be
characterized as its transfinite diameter, see e.g.  \cite{F} and
the point 110 in \cite{N}. Hence \begin{equation}\label{eqNP5}
C(f(Z_f))\ \leq C(f(U))\ \leq\ \varepsilon\cdot C(U)\ \leq\
\varepsilon\cdot C(D)\ .
\end{equation}
Thus,  by arbitrariness of $\varepsilon>0$ we obtain the right
equality in (\ref{eqNP3}).

Note that the function $J_f(z)=|f_z|^2-|f_{\bar z}|^2$ is continuous
and, consequently, the set $Z_f$ is closed in $D$, i.e. $Z_f$ and
$f(Z_f)$ are sigma-compacta in $D$ and $D_*=f(D)$, correspondingly.
Thus, the equality $C(f(Z_f))=0$ implies the equality $C(Z_f)=0$ by
the same characterization of logarithmic capacity because the
inverse mapping $f^{-1}$ is quasiconformal and hence it is H\"older
continuous on every compact set in the domain $D_*$, see e.g.
Theorem II.4.3 in \cite{LV}.
\end{proof}\ $\Box$

\bigskip

{\it Proof of Theorem 7.} By Lemma 1 $\mu$ is extended to a H\"older
continuous function $\mu_*:\mathbb C\to\mathbb C$. Hence also, for
every $k_*\in(k,1)$, there is a neighborhood $U$ of the closure of
the domain $D$ such that $|\mu(z)|<k_*$. Let $D_*$ be a connected
component of $U$ containing $\overline D$.

Now, there is a quasiconformal mapping $h:{{D_*}}\to{\mathbb{C}}$
a.e. satisfying the Beltrami equation \eqref{1} with the complex
coefficient $\mu^*:=\mu_*|_{D_*}$ in $D_*$, see e.g. Theorem V.B.3
in \cite{Alf}. Note that the mapping $h$ has the H\"older continuous
first partial derivatives in $D_*$ with the same order of the
H\"older continuity as $\mu$, see e.g. \cite{Iw} and also
\cite{IwDis}. Thus, by Lemma 2
\begin{equation}\label{eqNP6}
C(Z_h)\ =\ 0\ =\ C(h(Z_h))\ ,
\end{equation}
where $C(E)$ is the logarithmic capacity of $E\subseteq\mathbb C$,
$J_h$ is the Jacobian of $h$, and
\begin{equation}\label{eqNP7}
Z_h\ =\ \{\ z\in\mathbb C: J_h(z)\ =\ 0\ \}\ .
\end{equation}}
Moreover, setting $D^*=h(D)$ and
$\Gamma_{\xi}=h(\gamma_{h^{-1}(\xi)})$, $\xi\in\partial D^*$, we see
that $\{ \Gamma_{\xi}\}_{\xi\in\partial D^*}$ is a family of Jordan
arcs of class ${\cal{BS}}$ in $D^*$.

Next, setting
\begin{equation}\label{eqDERIVATIVE}
h_{\nu}(z)\ :=\ \lim_{t\to 0}\
\frac{h(z+t\cdot\nu)-h(z)}{t\cdot\nu}\ ,
\end{equation}
we see that the function $h_{\nu}(z)$ is continuous and hence
measurable and by (\ref{eqNP6}) also $h_{\nu}(z)\ne 0$ a.e. in $z\in
D_*$ with respect to the logarithmic capacity.

The logarithmic capacity of a set coincides with its transfinite
diameter, see e.g.  \cite{F} and the point 110 in \cite{N}. Hence
the mappings $h$ and $h^{-1}$ transform sets of logarithmic capacity
zero on $\partial D$ into sets of logarithmic capacity zero on
$\partial D^*$ and vice versa because these quasiconformal mappings
are continuous by H\"older on $\partial D$ and $\partial D^*$
correspondingly, see e.g. Theorem II.4.3 in \cite{LV}.

Further, the functions ${\cal{N}} := \nu\circ h^{-1}$ and $\varphi
:=\left(\Phi /h_{\nu}\right)\circ h^{-1}$ are mea\-su\-rab\-le with
respect to logarithmic capacity. Indeed, measurable sets with
respect to lo\-ga\-rith\-mic capacity are transformed under the
mappings $h$ and $h^{-1}$ into measurable sets with respect to
logarithmic capacity because such a set can be represented as the
union of a sigma-com\-pac\-tum and a set of logarithmic capacity
zero and compacta under continuous mappings are transformed into
compacta and compacta are measurable with respect to logarithmic
capacity.

Finally, by Theorem 6 there is a nonconstant single-valued analytic
function $F: D^*\to\mathbb C$ such that
\begin{equation}\label{NP29}
\lim_{w\to\xi}\ \frac{\partial F}{\partial {\cal{N}}}\ (w)\ =\
\varphi(\xi)
\end{equation}
holds along $\Gamma_{\xi}$ for a.e. $\xi\in\partial D^*$ with
respect to the logarithmic capacity.

It remains to note that $f:=F\circ h$ is the desired solution of
\eqref{1} because
$$
\frac{\partial f}{\partial\nu}\ =\ F^{\prime}\circ
h\cdot\frac{\partial h}{\partial\nu}\ =\ F^{\prime}\circ h\cdot\nu
h_{\nu}\ =\ \frac{\partial F}{\partial{\cal{N}}}\circ h\cdot
h_{\nu}\ ,
$$
i.e. (\ref{eqNP-LIMIT-R}) holds along $\gamma_{\zeta}$ for a.e.
$\zeta\in\partial D$ with respect to the logarithmic capacity.
$\Box$

\bigskip

{\bf Remark 5.} Again, we are able to say more in the case $\mathrm
{Re}\ n\cdot\overline{\nu}>0$ where $n=n(\zeta)$ is the unit
interior normal with a tangent to $\partial D$ at a point
$\zeta\in\partial D$. In view of (\ref{eqNP-LIMIT-R}), since the
limit $\Phi(\zeta)$ is finite, there is a finite limit $f(\zeta)$ of
$f(z)$ as $z\to\zeta$ in $D$ along the straight line passing through
the point $\zeta$ and being parallel to the vector $\nu(\zeta)$
because along this line, for $z$ and $z_0$ that are close enough to
$\zeta$,
\begin{equation}\label{eqDIFFERENCE-R} f(z)\ =\ f(z_0)\ -\ \int\limits_{0}\limits^{1}\
\frac{\partial f}{\partial \nu}\ (z_0+\tau (z-z_0))\ d\tau\
.\end{equation} Thus, at each point with the condition
(\ref{eqNP-LIMIT-R}), there is the directional derivative
\begin{equation}\label{eqPOSITIVE-R}
\frac{\partial f}{\partial \nu}\ (\zeta)\ :=\ \lim_{t\to 0}\
\frac{f(\zeta+t\cdot\nu)-f(\zeta)}{t}\ =\ \Phi(\zeta)\ .
\end{equation}

\medskip

In particular, we obtain on the basis of Theorem 7 and Remark 5 the
following statement on the Neumann problem.

\medskip

{\bf Corollary 5.} {\it\, Let $\mu:\mathbb D\to\mathbb{C}$ be in
$C^{\alpha}$, $\alpha\in(0,1)$, and $|\mu(z)|\leq k<1$, and let
$\Phi:\partial\mathbb D\to\mathbb C$ be measurable with respect to
the logarithmic capacity.

Then the Beltrami equation \eqref{1} has a regular solution $f:
\mathbb D\to\mathbb C$ of the class $C^{1+\alpha}$ such that, for
a.e. point $\zeta\in\partial\mathbb D$, there exist:

\bigskip

1) the finite radial limit
\begin{equation}\label{eqB-NP-LIMIT1}
f(\zeta)\ :=\ \lim\limits_{r\to 1}\ f(r\zeta)\end{equation}

2) the normal derivative
\begin{equation}\label{eqB-NP-NORMAL}
\frac{\partial f}{\partial n}\ (\zeta)\ :=\ \lim_{t\to 0}\
\frac{f(\zeta+t\cdot n)-f(\zeta)}{t}\ =\ \Phi(\zeta)
\end{equation}

3) the radial limit
\begin{equation}\label{eqB-NP-LIMIT2} \lim\limits_{r\to 1}\ \frac{\partial f}{\partial n}\ (r\zeta)\ =\
\frac{\partial f}{\partial n}\ (\zeta)\end{equation} where
$n=n(\zeta)$ denotes the unit interior normal to $\partial\mathbb D$
at the point $\zeta$.}

\medskip

{\bf Remark 6.} As it follows from Remark 4, Theorems 7 is valid for
multiply connected domains $D$ bounded by a finite collection of
mutually disjoint Jordan curves, however, then $f$ can be
multivalent in the spirit of the theory of multivalent analytic
functions.

\medskip

Namely, in  finitely connected domains $D$ in $\Bbb{C}$, we say that
a continuous discrete open mapping $f:B(z_0,\varepsilon_0)\to{\Bbb
C}$, where $B(z_0,\varepsilon_0)\subseteq D$, is a {\bf local
regular solution} of \eqref{1} if $f\in W_{\rm loc}^{1,1}$,
$J_f(z)\neq0$ and $f$ satisfies \eqref{1} a.e. in
$B(z_0,\varepsilon_0)$. Local regular solutions
$f_0:B(z_0,\varepsilon_0)\to{\Bbb C}$ and
$f_*:B(z_*,\varepsilon_*)\to{\Bbb C}$ of \eqref{1} are called the
extension of each to other if there is a finite chain of such
solutions $f_i:B(z_i,\varepsilon_i)\to\Bbb{C}$, $i=1,\ldots,m$, that
$f_1=f_0$, $f_m=f_*$ and $f_i(z)\equiv f_{i+1}(z)$ for $z\in
E_i:=B(z_i,\varepsilon_i)\cap
B(z_{i+1},\varepsilon_{i+1})\neq\emptyset$, $i=1,\ldots,m-1$. A
collection of local regular solutions
$f_j:B(z_j,\varepsilon_j)\to{\Bbb C}$, $j\in J$, will be called a
{\bf multi-valued solution} of the Beltrami equation \eqref{1} in
$D$ if the disks $B(z_j,\varepsilon_j)$ cover the whole domain $D$
and $f_j$ are extensions of each to other through the collection and
this collection is maximal by inclusion, cf. e.g. \cite{BGR} and
\cite{KPRS}.

\bigskip

\section{On the Reimann problem for analytic functions}

\medskip

Recall that the classical setting of the {\bf Riemann problem} in a
smooth Jordan domain $D$ of the complex plane $\mathbb{C}$ was on
finding analytic functions $f^+: D\to\mathbb C$ and $f^-:\mathbb
C\setminus \overline{D}\to\mathbb C$ that admit continuous
extensions to $\partial D$ and satisfy the boundary condition
\begin{equation}\label{eqRIEMANN} f^+(\zeta)\ =\
A(\zeta)\cdot f^-(\zeta)\ +\ B(\zeta) \quad\quad\quad \forall\
\zeta\in\partial D
\end{equation}
with prescribed H\"older continuous functions $A: \partial
D\to\mathbb C$ and $B: \partial D\to\mathbb C$.

\medskip

Recall also that the {\bf Riemann problem with shift} in $D$ was on
finding analytic functions $f^+: D\to\mathbb C$ and $f^-:\mathbb
C\setminus \overline{D}\to\mathbb C$ satisfying the condition
\begin{equation}\label{eqSHIFT} f^+(\alpha(\zeta))\ =\ A(\zeta)\cdot
f^-(\zeta)\ +\ B(\zeta) \quad\quad\quad \forall\ \zeta\in\partial D
\end{equation}
where $\alpha :\partial D\to\partial D$ was a one-to-one sense
preserving correspondence having the non-vanishing H\"older
continuous derivative with respect to the natural parameter on
$\partial D$. The function $\alpha$ is called a {\bf shift
function}. The special case $A\equiv 1$ gives the so--called {\bf
jump problem} and $B\equiv 0$ gives the {\bf problem on gluing} of
analytic functions.

\bigskip

{\bf Theorem 8.} {\it\, Let $D$ be a domain in $\overline{\mathbb
C}$ whose boundary consists of a finite number of mutually disjoint
Jordan curves, $A: \partial D\to\mathbb C$ and $B:
\partial D\to\mathbb C$ be functions that are measurable with respect to the
logarithmic capacity. Suppose that
$\{\gamma^+_{\zeta}\}_{\zeta\in\partial D}$ and
$\{\gamma^-_{\zeta}\}_{\zeta\in\partial D}$ are families of Jordan
arcs of class ${\cal{BS}}$ in ${D}$ and $\mathbb
C\setminus\overline{ D}$, correspondingly.

\medskip

Then there exist single-valued analytic functions $f^+: D\to\mathbb
C$ and $f^-:\overline{\mathbb C}\setminus\overline{D}\to\mathbb C$
that satisfy (\ref{eqRIEMANN}) for a.e. $\zeta\in\partial D$ with
respect to the logarithmic capacity where $f^+(\zeta)$ and
$f^-(\zeta)$ are limits of $f^+(z)$ and $f^-(z)$ az $z\to\zeta$
along $\gamma^+_{\zeta}$ and $ \gamma^-_{\zeta}$, correspondingly.

\medskip

Furthermore, the space of all such couples $(f^+,f^-)$ has the
infinite dimension for every couple $(A, B)$ and any collections
$\gamma^+_{\zeta}$ and $ \gamma^-_{\zeta}$, $\zeta\in\partial D$.}

\medskip

Theorem 8 is a special case of the following lemma on the
generalized Riemann problem with shift that can be useful for other
goals, too.

\medskip

{\bf Lemma 3.} {\it\, Under the hypotheses of Theorem 8, let in
addition $\alpha : \partial D\to\partial D$ be a homeomorphism
keeping components of $\partial D$ such that $\alpha$ and
$\alpha^{-1}$ have the $(N)-$property of Lusin with respect to the
logarithmic capacity.

\medskip

Then there exist single-valued analytic functions $f^+: D\to\mathbb
C$ and $f^-:\overline{\mathbb C}\setminus\overline{D}\to\mathbb C$
that satisfy (\ref{eqSHIFT}) for a.e. $\zeta\in\partial D$ with
respect to the logarithmic capacity where $f^+(\zeta)$ and
$f^-(\zeta)$ are limits of $f^+(z)$ and $f^-(z)$ az $z\to\zeta$
along $\gamma^+_{\zeta}$ and $ \gamma^-_{\zeta}$, correspondingly.

\medskip

Furthermore, the space of all such couples $(f^+,f^-)$ has the
infinite dimension for every couple $(A, B)$ and any collections
$\gamma^+_{\zeta}$ and $ \gamma^-_{\zeta}$, $\zeta\in\partial D$.}

\bigskip

\begin{proof} For short, we write here "$C-$measurable" instead of
the expression "measurable with respect to the logarithmic
capacity".

First, let $D$ be bounded and let $g^-: \partial D\to\mathbb C$ be a
$C-$measurable function. Then the function
\begin{equation}\label{eqCONNECTION}g^+\ :=\ \{ A\cdot
g^- + B\}\circ \alpha^{-1}\end{equation} is also $C-$measurable.
Indeed, $E:=\{ A\cdot g^- + B\}^{-1}(\Omega)$ is a subset of
$\partial D$ that is $C-$measurable for every open set
$\Omega\subseteq\mathbb C$ because the function $A\cdot g^- + B$ is
$C-$measurable by the hypotheses. Hence the set $E$ is the union of
a sigma-compact set and a set of the logarithmic capacity zero, see
e.g. Section 2 in \cite{ER}. However, continuous mappings transform
compact sets into compact sets and, thus, by $(N)-$property
$\alpha(E)=\alpha\circ\{ A\cdot g^- + B\}
^{-1}(\Omega)=(g^+)^{-1}(\Omega)$ is a $C-$measurable set, i.e. the
function $g^+$ is really $C-$measurable.

Then by Theorem 1 there is a single-valued analytic function $f^+:
D\to\mathbb C$ such that
\begin{equation}\label{eqBS-R2} \lim\limits_{z\to\xi}\ f^+(z)\ =\ g^+(\xi)
\end{equation}
along $\gamma^+_{\xi}$ for a.e. $\xi\in\partial D$ with respect to
the logarithmic capacity. Note that $g^+(\alpha(\zeta))$ is
determined by the given limit for a.e. $\zeta\in\partial D$ because
$\alpha^{-1}$ also has the $(N)-$property of Lusin.

Note that $\overline{\mathbb C}\setminus\overline{D}$ consists of a
finite number of (simply connected) Jordan domains $D_0, D_1, \ldots
, D_m$ in the extended complex plane $\overline{\mathbb C}=\mathbb
C\cup \{ \infty\}$. Let $\infty\in D_0$. Then again by Theorem 1
there exist single-valued analytic functions $f_l^-: D_l\to\mathbb
C$, $l=1, \ldots , m,$ such that
\begin{equation}\label{eqBS-R1} \lim\limits_{z\to\zeta}\ f_l^-(z)\ =\
g_l^-(\zeta)\ ,\quad\quad\quad g_l^-:=g^-|_{\partial D_l}\ ,
\end{equation}
along $\gamma^-_{\zeta}$ for a.e. $\zeta\in\partial D_l$ with
respect to the logarithmic capacity.

Now, let $S$ be a circle that contains $D$ and let $j$ be the
inversion of $\overline{\mathbb C}$ with respect to $S$. Set
$$D_*=j(D_0),\quad g_*=\overline {g_0\circ j},\quad  g_0^-:=g^-|_{\partial
D_0},\quad \gamma^*_{\xi}=j\left(\gamma^-_{j(\xi)}\right),\quad
\xi\in\partial D_*\ .$$ Then by Theorem 1 there is a single-valued
analytic function $f_*: D_*\to\mathbb C$ such that
\begin{equation}\label{eqBS-R} \lim\limits_{w\to\xi}\ f_*(w)\ =\ g_*(\xi)
\end{equation}
along $\gamma^*_{\xi}$ for a.e. $\xi\in\partial D_*$ with respect to
the logarithmic capacity. Note that $f_0^-:=\overline {g_*\circ j}$
is a single-valued analytic function in $D_0$ and by construction
\begin{equation}\label{eqBS0-R} \lim\limits_{z\to\zeta}\ f_0^-(z)\ =\
g_0^-(\zeta)\ ,\quad\quad\quad g_0^-:=g^-|_{\partial D_0}\ ,
\end{equation} along $\gamma^-_{\zeta}$ for a.e. $\zeta\in\partial D_0$ with
respect to the logarithmic capacity.

Thus, the functions $f^-_l$, $l=0, 1, \ldots , m,$ form an analytic
function $f^-: \overline{\mathbb C}\setminus\overline{D}\to\mathbb
C$ satisfying (\ref{eqSHIFT}) for a.e. $\zeta\in\partial D$ with
respect to the logarithmic capacity.

The space of all such couples $(f^+,f^-)$ has the infinite dimension
for every couple $(A, B)$ and any collections $\gamma^+_{\zeta}$ and
$ \gamma^-_{\zeta}$, $\zeta\in\partial D$, in view of the above
construction because of the space of all measurable functions $g^-:
\partial D\to\mathbb C$ has the infinite dimension.

The case of unbounded $D$ is reduced to the case of bounded $D$
through the complex conjugation and the inversion of
$\overline{\mathbb C}$ with respect to a circle $S$ in some of the
components of $\overline{\mathbb C}\setminus\overline{D}$ arguing as
above.
\end{proof}

\bigskip

{\bf Remark 7.} Some investigations were devoted also to the
nonlinear Riemann problems with boundary conditions of the form
\begin{equation}\label{eqNONLINEAR} \Phi(\,\zeta,\, f^+(\zeta),\, f^-(\zeta)\, )\ =\ 0 \quad\quad\quad \forall\
\zeta\in\partial D\ .
\end{equation}
It is natural as above to weaken such conditions to the following
\begin{equation}\label{eqNONLINEAR} \Phi(\,\zeta,\, f^+(\zeta),\, f^-(\zeta)\, )\ =\ 0 \quad\quad\quad \mbox{for a.e.}\quad
\zeta\in\partial D\ .
\end{equation}
It is easy to see that the proposed approach makes possible  also to
reduce such problems to the algebraic and measurable solvability of
the relations
\begin{equation}\label{eqNONLINEAR} \Phi(\,\zeta,\, v,\, w\, )\ =\ 0
\end{equation}
with respect to complex-valued functions $v(\zeta)$ and $w(\zeta)$,
cf. e.g. \cite{Gr}.

\bigskip

{\bf Example 1.} For instance, correspondingly to the scheme given
above, special nonlinear problems of the form
\begin{equation}\label{eqNONLINEAR} f^+(\zeta)\ =\ \varphi(\,\zeta,\,
 f^-(\zeta)\, ) \quad\quad\quad \mbox{for
a.e.}\quad \zeta\in\partial D
\end{equation}
are always solved if the function $\varphi : \partial D\times\mathbb
C\to\mathbb C$ satisfies the {\bf Caratheodory conditions} with
respect to the logarithmic capacity: $\varphi(\zeta, w)$ is
continuous in the variable $w\in\mathbb C$ for a.e.
$\zeta\in\partial D$ with respect to the logarithmic capacity and it
is $C-$measurable in the variable $\zeta\in \partial D$ for all
$w\in\mathbb C$.

\medskip

Furthermore, the spaces of solutions of such problems always have
the infinite dimension. Indeed, by the Egorov theorem, see e.g.
Theorem 2.3.7 in \cite{Fe}, see also Section 17.1 in \cite{KZPS},
the function $\varphi(\zeta,\psi(\zeta))$ is $C-$measurable in
$\zeta\in\partial D$ for every $C-$measurable function
$\psi:\partial D\to\mathbb C$ if the function $\varphi$ satisfies
the {Caratheodory conditions}, and the space of all $C-$measurable
functions $\psi:\partial D\to\mathbb C$ has the infinite dimension,
see e.g. Remark 1.\bigskip

\section{On the Riemann problem for the Beltrami equations}

{\bf Theorem 9.} {\it\, Let $D$ be a domain in $\overline{\mathbb
C}$ whose boundary consists of a finite number of mutually disjoint
Jordan curves, $\mu:\mathbb{C}\to\mathbb{C}$ be a measurable (by
Lebesgue) function with $||\mu||_{\infty}<1$, $A: \partial
D\to\mathbb C$ and $B:
\partial D\to\mathbb C$ be functions that are measurable with respect to the
logarithmic capacity. Suppose that
$\{\gamma^+_{\zeta}\}_{\zeta\in\partial D}$ and
$\{\gamma^-_{\zeta}\}_{\zeta\in\partial D}$ are families of Jordan
arcs of class ${\cal{BS}}$ in ${D}$ and $\mathbb
C\setminus\overline{ D}$, correspondingly.

\medskip

Then the Beltrami equation \eqref{1} has regular solutions $f^+:
D\to\mathbb C$ and $f^-:\overline{\mathbb
C}\setminus\overline{D}\to\mathbb C$ that satisfy (\ref{eqRIEMANN})
for a.e. $\zeta\in\partial D$ with respect to the logarithmic
capacity where $f^+(\zeta)$ and $f^-(\zeta)$ are limits of $f^+(z)$
and $f^-(z)$ az $z\to\zeta$ along $\gamma^+_{\zeta}$ and $
\gamma^-_{\zeta}$, correspondingly.

\medskip

Furthermore, the space of all such couples $(f^+,f^-)$ has the
infinite dimension for every couple $(A, B)$ and any collections
$\gamma^+_{\zeta}$ and $ \gamma^-_{\zeta}$, $\zeta\in\partial D$.}

\medskip

Theorem 9 is a special case of the following lemma on the Riemann
problem with shift for the Beltrami equation.

\medskip

{\bf Lemma 4.} {\it\, Under the hypotheses of Theorem 9, let in
addition $\alpha : \partial D\to\partial D$ be a homeomorphism
keeping components of $\partial D$ such that $\alpha$ and
$\alpha^{-1}$ have the $(N)-$property of Lusin with respect to the
logarithmic capacity.

\medskip

Then the Beltrami equation \eqref{1} has regular solutions $f^+:
D\to\mathbb C$ and $f^-:\overline{\mathbb
C}\setminus\overline{D}\to\mathbb C$ that satisfy (\ref{eqSHIFT})
for a.e. $\zeta\in\partial D$ with respect to the logarithmic
capacity where $f^+(\zeta)$ and $f^-(\zeta)$ are limits of $f^+(z)$
and $f^-(z)$ az $z\to\zeta$ along $\gamma^+_{\zeta}$ and $
\gamma^-_{\zeta}$, correspondingly.

\medskip

Furthermore, the space of all such couples $(f^+,f^-)$ has the
infinite dimension for every couple $(A, B)$ and any collections
$\gamma^+_{\zeta}$ and $ \gamma^-_{\zeta}$, $\zeta\in\partial D$.}

\bigskip

\begin{proof} First of all, there is a quasiconformal mapping
$h:\mathbb{C}\to\mathbb{C}$ that is a regular homeomorphic solution
of the Beltrami equation \eqref{1} with the given complex
coefficient $\mu$, see e.g. Theorem V.B.3 in \cite{Alf}. Set
$D^*=h(D)$ and $\alpha_*=h\circ\alpha\circ h^{-1}$. Then $\alpha_* :
\partial D_*\to\partial D_*$ is a homeomorphism.

Recall that the logarithmic capacity of a set coincides with its
transfinite diameter, see \cite{F} and the point 110 in \cite{N}.
Hence the mappings $h$ and $h^{-1}$ transform sets of logarithmic
capacity zero on $\partial D$ into sets of logarithmic capacity zero
on $\partial D^*$ and vice versa because quasiconformal mappings are
continuous by H\"older on $\partial D$ and $\partial D^*$
correspondingly, see e.g. Theorem II.4.3 in \cite{LV}. Thus,
$\alpha_*$ and $\alpha^{-1}_*$ have the $(N)-$property with respect
to the logarithmic capacity.

Moreover, arguing similarly to the proof of Theorem 7, it is easy to
show that the functions $A_*:=A\circ h^{-1}$ and $B_*:=B\circ
h^{-1}$ are measurable with respect to the logarithmic capacity.
Setting $\Gamma^+_{\xi}=h(\gamma^+_{h^{-1}(\xi)})$ and
$\Gamma^-_{\xi}=h(\gamma^-_{h^{-1}(\xi)})$, $\xi\in\partial D^*$, we
see also that $\{\Gamma^+_{\xi}\}_{\xi\in\partial D^*}$ and
$\{\Gamma^-_{\xi}\}_{\xi\in\partial D^*}$ are families of Jordan
arcs of class ${\cal{BS}}$ in $D^*$  and $\mathbb
C\setminus\overline{ D^*}$, correspondingly.

Hence by Lemma 3 there exist single-valued analytic functions
$f_*^+: D_*\to\mathbb C$ and $f_*^-:\overline{\mathbb
C}\setminus\overline{D_*}\to\mathbb C$ such that
\begin{equation}\label{eqSHIFT-B} f_*^+(\alpha_*(\xi))\ =\ A_*(\xi)\cdot
f_*^-(\xi)\ +\ B_*(\xi) \quad\quad\quad \forall\ \xi\in\partial D_*
\end{equation}
for a.e. $\xi\in\partial D_*$ with respect to the logarithmic
capacity where $f_*^+(\xi)$ and $f_*^-(\xi)$ are limits of
$f_*^+(w)$ and $f_*^-(w)$ as $w\to\xi$ along $\Gamma^+_{\xi}$ and $
\Gamma^-_{\xi}$. Furthermore, the space of all such couples
$(f_*^+,f_*^-)$ has the infinite dimension.

Finally, $f^+:=f^+_*\circ h$ and  $f^-:=f^-_*\circ h$ are the
desired regular solutions of the Beltrami equation \eqref{1} in
$D^*$  and $\mathbb C\setminus\overline{ D^*}$, correspondingly.
\end{proof} $\Box$

\bigskip

{\bf Remark 7.} As it follows from the construction of solutions in
the given proof, Remark 6 on the nonlinear Riemann problems is valid
for the Beltrami equations, see also Example 1 in the last section.

\section{On mixed boundary value problems}

For the demonstration of potentialities latent in our approach, we
consider here some nonlinear boundary value problems.

\medskip

{\bf Theorem 10.} {\it\, Let $D$ be a domain in $\overline{\mathbb
C}$ whose boundary consists of a finite number of mutually disjoint
Jordan curves, $\mu:\mathbb{C}\to\mathbb{C}$ be measurable (by
Lebesgue), $\mu\in C^{\alpha}(\mathbb C\setminus\overline{ D})$ and
$||\mu||_{\infty}<1$, $\varphi :
\partial D\times\mathbb C\to\mathbb C$ satisfy the Caratheodory
conditions and $\nu:\partial D\to\mathbb C$, $|\nu (\zeta)|\equiv
1$, be measurable with respect to the logarithmic capacity.

\medskip

Suppose that $\{\gamma^+_{\zeta}\}_{\zeta\in\partial D}$ and
$\{\gamma^-_{\zeta}\}_{\zeta\in\partial D}$ are families of Jordan
arcs of class ${\cal{BS}}$ in ${D}$ and $\mathbb
C\setminus\overline{ D}$. Then the Beltrami equation \eqref{1} has
regular solutions $f^+: D\to\mathbb C$ and $f^-:\overline{\mathbb
C}\setminus\overline{D}\to\mathbb C$ of class $C^{1+\alpha}(\mathbb
C\setminus\overline{ D})$ such that
\begin{equation}\label{eqMIXED} f^+(\zeta)\ =\ \varphi\left(\,\zeta,\,
\left[\frac{\partial f}{\partial\nu}\right]^- (\zeta)\, \right)
\end{equation}
for a.e. $\zeta\in\partial D$ with respect to the logarithmic
capacity where $f^+(\zeta)$ and $\left[\frac{\partial
f}{\partial\nu}\right]^- (\zeta)$ are limits of the functions
$f^+(z)$ and $\frac{\partial f^-}{\partial\, \nu\,}\ (z)$ as
$z\to\zeta$ along $\gamma^+_{\zeta}$ and $ \gamma^-_{\zeta}$,
correspondingly.

\medskip

Furthermore, the space of all such couples $(f^+,f^-)$ has the
infinite dimension for any such prescribed functions $\mu$,
$\varphi$, $\nu$ and collections $\gamma^+_{\zeta}$ and $
\gamma^-_{\zeta}$, $\zeta\in\partial D$.}

\medskip

Theorem 10 is a special case of the following lemma on the nonlinear
Riemann problem with shift.

\medskip

{\bf Lemma 5.} {\it\, Under the hypotheses of Theorem 10, let in
addition $\beta : \partial D\to\partial D$ be a homeomorphism
keeping components of $\partial D$ such that $\beta$ and
$\beta^{-1}$ have the $(N)-$property of Lusin with respect to the
logarithmic capacity.

\medskip

Then the Beltrami equation \eqref{1} has regular solutions $f^+:
D\to\mathbb C$ and $f^-:\overline{\mathbb
C}\setminus\overline{D}\to\mathbb C$ of class $C^{1+\alpha}(\mathbb
C\setminus\overline{ D})$ such that
\begin{equation}\label{eqMIXED_SHIFT} f^+(\beta(\zeta))\ =\
\varphi\left(\,\zeta,\, \left[\frac{\partial
f}{\partial\nu}\right]^- (\zeta)\, \right)
\end{equation}
for a.e. $\zeta\in\partial D$ with respect to the logarithmic
capacity where $f^+(\zeta)$ and $\left[\frac{\partial
f}{\partial\nu}\right]^- (\zeta)$ are limits of the functions
$f^+(z)$ and $\frac{\partial f^-}{\partial\, \nu\,}\ (z)$ as
$z\to\zeta$ along $\gamma^+_{\zeta}$ and $ \gamma^-_{\zeta}$,
correspondingly.

\medskip

Furthermore, the space of all such couples $(f^+,f^-)$ has the
infinite dimension for any such prescribed $\mu$, $\varphi$, $\nu$,
$\beta$ and collections $\{\gamma^+_{\zeta}\}_{\zeta\in\partial D}$
and $\{\gamma^-_{\zeta}\}_{\zeta\in\partial D}$.}

\bigskip

\begin{proof} Indeed, let $\psi:\partial D\to\mathbb C$ be a
function that is measurable with respect to the logarithmic
capacity. Then by Theorem 7 the Beltrami equation \eqref{1} has a
regular solution $f^-:\overline{\mathbb
C}\setminus\overline{D}\to\mathbb C$ of class $C^{1+\alpha}$ such
that
\begin{equation}\label{eqNP-LIMIT-MIXED} \lim\limits_{z\to\zeta}\
\frac{\partial f^-}{\partial\, \nu\,}\ (z)\ =\
\psi(\zeta)\end{equation} along $\gamma^-_{\zeta}$ for a.e.
$\zeta\in\partial D$ with respect to the logarithmic capacity.

Now, the function $\ \Psi(\zeta)\ :=\ \varphi(\zeta,\psi(\zeta))\ $
is measurable with respect to the logarithmic capacity, see Example
1 after Remark 6. Then the function $\Phi = \Psi\circ\alpha^{-1}$ is
also measurable with respect to the logarithmic capacity because the
homeomorphism $\alpha$ has the $(N)-$property, cf. arguments in the
proof of Lemma 4.

Thus, by Theorem 2 the Beltrami equation \eqref{1} has a regular
solution $f^+: D\to\mathbb C$ such that $f^+(z)\to\Psi(\zeta)$ as
$z\to\zeta$ along $\gamma_{\zeta}$ for a.e. $\zeta\in\partial D$
with respect to the logarithmic capacity. Finally, $f^+$ and $f^-$
are the desired functions because $\alpha^{-1}$ also has the
$(N)-$property.

It remains also to note that the space of all such couples
$(f^+,f^-)$ has the infinite dimension because the space of all
functions $\psi:\partial D\to\mathbb C$ which are measurable with
respect to the logarithmic capacity has the infinite dimension, see
Remark 1.
\end{proof} $\Box$

\section{On applications to equations of the divergence type}

The partial differential equations in the divergence form below take
a significant part in many problems of mathematical physics, in
particular, in anisotropic inhomogeneous media, see e.g. \cite{AIM},
\cite{IwHar} and \cite{IwSb}.

\bigskip

In this connection, note that if $f=u+i\cdot v$ is regular solution
of the Beltrami equation \eqref{1}, then the function $u$ is a
continuous generalized solution of the divergence-type equation
\begin{equation}\label{homo}
{\rm div}\, A(z)\nabla\,u=0\ ,
\end{equation}
called {\bf A-harmonic function}, i.e. $u\in C\cap W^{1,1}$ and
\begin{equation}\label{weak}
\int\limits_D \langle A(z)\nabla
u,\nabla\varphi\rangle=0\,\,\,\,\,\,\,\,\,\,\,\,\forall\ \varphi\in
C_0^\infty(D)\ ,
\end{equation}
where $A(z)$ is the matrix function:
\begin{equation}\label{matrix}
A=\left(\begin{array}{ccc} {|1-\mu|^2\over 1-|\mu|^2}  & {-2{\rm Im}\,\mu\over 1-|\mu|^2} \\
                            {-2{\rm Im}\,\mu\over 1-|\mu|^2}          & {|1+\mu|^2\over 1-|\mu|^2}  \end{array}\right).
\end{equation}
Moreover, we have for the {\bf stream function} $v={\rm Im}\, f$ the
following relation
\begin{equation}\label{con}
\nabla\, v=J\, A(z)\nabla\,u\ ,
\end{equation}
where
\begin{equation}\label{H}
J=\left(\begin{array}{rr} 0 &-1 \\
        1 &0\end{array}\right)\ .
\end{equation}
As we see in (\ref{matrix}), the matrix $A(z)$ is symmetric, ${\rm
det}\,A(z)\equiv 1$ and its entries $a_{ij}=a_{ij}(z)$ are dominated
by the quantity
\begin{equation}\label{dilatation}
K_{\mu}(z)\ =\ \frac{1\ +\ |\mu(z)|}{1\ -\ |\mu(z)|}\ ,
\end{equation}
and, thus, they are bounded if the Beltrami equation \eqref{1} is
not degenerate.

Vice verse, uniformly elliptic equations (\ref{homo}) with symmetric
$A(z)$ and ${\rm det}\,A(z)\equiv 1$ just correspond to
nondegenerate Beltrami equations \eqref{1}. We call such matrix
functions $A(z)$ of {\bf class ${\cal{B}}$}.

Recall that the equation (\ref{homo}) is called {\bf uniformly
elliptic} if $a_{ij}\in L^{\infty}$ and there is $\varepsilon > 0$
such that $\langle A(z)\eta,\eta\rangle\geq\varepsilon |\eta|^2$ for
all $\eta\in\mathbb R^2$. A matrix function $A(z)$ with the latter
property is called {\bf uniformly positive definite}.

\bigskip

The following statement on a potential boundary behavior of
A-harmonic functions of the Dirichlet type is a direct consequence
of Theorem 2.

\medskip

{\bf Corollary 6.} {\it  Let $D$ be a bounded domain in $\mathbb C$
whose boundary consists of a finite number of mutually disjoint
Jordan curves, $A(z)$, $z\in D$, be a matrix function of class
${\cal{B}}$ and let a function $\varphi:\partial D\to\mathbb R$ be
measurable with respect to the logarithmic capacity.

Suppose that $\{ \gamma_{\zeta}\}_{\zeta\in\partial D}$ is a family
of Jordan arcs of class ${\cal{BS}}$ in ${D}$. Then there exist
A-harmonic functions $u: D\to\mathbb R$ such that
\begin{equation}\label{eqBS-A} \lim\limits_{z\to\zeta}\ u(z)\ =\ \varphi(\zeta)
\end{equation}
along $\gamma_{\zeta}$ for a.e. $\zeta\in\partial D$ with respect to
the logarithmic capacity.

Furthermore, the space of all such A-harmonic functions $u$ has the
infinite dimension for any such prescribed $A$, $\varphi$ and $\{
\gamma_{\zeta}\}_{\zeta\in\partial D}$.}

\bigskip

The next conclusion in the particular case of the Poincare problem
on directional derivatives follows directly from Theorem 7.

\medskip

{\bf Corollary 7.} {\it\, Let $D$ be a Jordan domain in $\mathbb C$,
$A(z)$, $z\in D$, be a matrix function of class ${\cal{B}}\cap
C^{\alpha}$, $\alpha\in(0,1)$, and let $\nu:\partial D\to\mathbb C$,
$|\nu (\zeta)|\equiv 1$, and $\varphi:\partial D\to\mathbb R$ be
measurable with respect to the logarithmic capacity.

\medskip

Suppose that $\{ \gamma_{\zeta}\}_{\zeta\in\partial D}$ is a family
of Jordan arcs of class ${\cal{BS}}$ in ${D}$. Then there exist
A-harmonic functions $u: D\to\mathbb R$ of the class $C^{1+\alpha}$
such that
\begin{equation}\label{eqNP-LIMIT-R-A} \lim\limits_{z\to\zeta}\ \frac{\partial u}{\partial \nu}\ (z)\ =\
\varphi(\zeta)\end{equation} along $\gamma_{\zeta}$ for a.e.
$\zeta\in\partial D$ with respect to the logarithmic capacity.

Furthermore, the space of all such A-harmonic functions $u$ has the
infinite dimension for any such prescribed $A$, $\varphi$, $\nu$ and
$\{ \gamma_{\zeta}\}_{\zeta\in\partial D}$.}

\bigskip

Now, the following statement concerning to the Neumann problem for
A-harmonic functions is a special significant case of Corollary 7.

\medskip

{\bf Corollary 8.} {\it\, Let $A(z)$, $z\in\mathbb D$, be a matrix
function of class ${\cal{B}}\cap C^{\alpha}$, $\alpha\in(0,1)$, and
let $\varphi:\partial\mathbb D\to\mathbb R$ be measurable with
respect to the logarithmic capacity.

Then there exist A-harmonic functions $u:\mathbb D\to\mathbb R$ of
the class $C^{1+\alpha}$ such that, for a.e. point
$\zeta\in\partial\mathbb D$ with respect to the logarithmic
capacity, there exist:

\bigskip

1) the finite radial limit
\begin{equation}\label{eqB-NP-LIMIT1-A}
u(\zeta)\ :=\ \lim\limits_{r\to 1}\ u(r\zeta)\end{equation}

2) the normal derivative
\begin{equation}\label{eqB-NP-NORMAL-A}
\frac{\partial u}{\partial n}\ (\zeta)\ :=\ \lim_{t\to 0}\
\frac{u(\zeta+t\cdot n)-u(\zeta)}{t}\ =\ \varphi(\zeta)
\end{equation}

3) the radial limit
\begin{equation}\label{eqB-NP-LIMIT2-A} \lim\limits_{r\to 1}\ \frac{\partial u}{\partial n}\ (r\zeta)\ =\
\frac{\partial u}{\partial n}\ (\zeta)\end{equation} where
$n=n(\zeta)$ denotes the unit interior normal to $\partial\mathbb D$
at the point $\zeta$.

Furthermore, the space of all such A-harmonic functions $u$ has the
infinite dimension for any such prescribed $A$ and $\varphi$.}

\bigskip

The following result on the Riemann problem with shift for
A-harmonic functions is not a direct consequence of Lemma 4 but of
its proof where we are able to choose $f^-(\zeta)$ such that ${\rm
Im} f^-(\zeta)=0$ for a.e. $\zeta\in\partial D$ with respect to the
logarithmic capacity.

\medskip

{\bf Corollary 9.} {\it\, Let $D$ be a domain in ${\mathbb C}$ whose
boundary consists of a finite number of mutually disjoint Jordan
curves, $A(z)$, $z\in D$, be a matrix function of class ${\cal{B}}$,
$B: \partial D\to\mathbb R$ and $C:
\partial D\to\mathbb R$ be functions that are measurable with respect to the
logarithmic capacity and let $\alpha : \partial D\to\partial D$ be a
homeomorphism keeping components of $\partial D$ such that $\alpha$
and $\alpha^{-1}$ have the $(N)-$property of Lusin with respect to
the logarithmic capacity.

\medskip

Suppose that $\{\gamma^+_{\zeta}\}_{\zeta\in\partial D}$ and
$\{\gamma^-_{\zeta}\}_{\zeta\in\partial D}$ are families of Jordan
arcs of class ${\cal{BS}}$ in ${D}$ and $\mathbb
C\setminus\overline{ D}$, correspondingly. Then there exist
A-harmonic functions $u^+: D\to\mathbb R$ and $u^-:\overline{\mathbb
C}\setminus\overline{D}\to\mathbb R$ such that
\begin{equation}\label{eqSHIFT-A} u^+(\alpha(\zeta))\ =\ B(\zeta)\cdot
u^-(\zeta)\ +\ C(\zeta)
\end{equation}
for a.e. $\zeta\in\partial D$ with respect to the logarithmic
capacity where $u^+(\zeta)$ and $u^-(\zeta)$ are limits of $u^+(z)$
and $u^-(z)$ az $z\to\zeta$ along $\gamma^+_{\zeta}$ and $
\gamma^-_{\zeta}$, correspondingly.

\medskip

Furthermore, the space of all such couples $(u^+,u^-)$ has the
infinite dimension for any such prescribed $A$, $B$, $C$, $\alpha$
and collections $\{\gamma^+_{\zeta}\}_{\zeta\in\partial D}$ and
$\{\gamma^-_{\zeta}\}_{\zeta\in\partial D}$.}

\bigskip

Similarly, the following result on the mixed nonlinear Riemann
problem with shifts for A-harmonic functions is not a direct
consequence of Lemma 5 but of its proof where we may choose
$f^-(\zeta)$ such that ${\rm Im} \left[\frac{\partial
f}{\partial\nu}\right]^- (\zeta)=0$ for a.e. $\zeta\in\partial D$
with respect to the logarithmic capacity.

\medskip

{\bf Corollary 10.} {\it\, Let $D$ be a domain in ${\mathbb C}$
whose boundary consists of a finite number of mutually disjoint
Jordan curves, $A(z)$, $z\in \mathbb C$, be a matrix function of
class  ${\cal{B}}\cap C^{\alpha}(\mathbb C\setminus\overline{ D})$,
$\alpha\in(0,1)$,  $\nu:\partial D\to\mathbb C$, $|\nu
(\zeta)|\equiv 1$, be a measurable function, $\beta :
\partial D\to\partial D$ be a homeomorphism such that $\beta$ and
$\beta^{-1}$ have the $(N)-$property of Lusin and $\varphi :
\partial D\times\mathbb R\to\mathbb R$ satisfy the Caratheodory
conditions with respect to the logarithmic capacity.

\medskip

Suppose that $\{\gamma^+_{\zeta}\}_{\zeta\in\partial D}$ and
$\{\gamma^-_{\zeta}\}_{\zeta\in\partial D}$ are families of Jordan
arcs of class ${\cal{BS}}$ in ${D}$ and $\mathbb
C\setminus\overline{ D}$, correspondingly. Then there exist
A-harmonic functions $u^+: D\to\mathbb R$ and $u^-:{\mathbb
C}\setminus\overline{D}\to\mathbb R$ of class $C^{1+\alpha}(\mathbb
C\setminus\overline{ D})$ such that
\begin{equation}\label{eqMIXED-A} u^+(\beta(\zeta))\ =\ \varphi\left(\,\zeta,\,
\left[\frac{\partial u}{\partial\nu}\right]^- (\zeta)\, \right)
\end{equation}
for a.e. $\zeta\in\partial D$ with respect to the logarithmic
capacity where $u^+(\zeta)$ and $\left[\frac{\partial
u}{\partial\nu}\right]^- (\zeta)$ are limits of the functions
$u^+(z)$ and $\frac{\partial u^-}{\partial\, \nu\,}\ (z)$ as
$z\to\zeta$ along $\gamma^+_{\zeta}$ and $ \gamma^-_{\zeta}$,
correspondingly.

\medskip

Furthermore, the space of all such couples $(u^+,u^-)$ has the
infinite dimension for any such prescribed $A$, $\nu$, $\beta$,
$\varphi$ and collections $\gamma^+_{\zeta}$ and $
\gamma^-_{\zeta}$, $\zeta\in\partial D$.}

\bigskip

In particular, we are able to obtain from the last corollary
solutions of the problem on gluing of the Dirichlet problem in the
unit disk $\mathbb D$ and the Neumann problem outside of $\mathbb D$
in the class of A-harmonic functions.

\medskip
\noindent {\bf Vladimir Gutlyanskii,
Vladimir Ryazanov,\\
Artem Yefimushkin}\\
Institute of Applied Mathematics and Mechanics,\\
National Academy of Sciences of Ukraine,\\
1 Dobrovolskogo Str., Slavyansk, 84100, UKRAINE,\\
vladimirgut@mail.ru, vl.ryazanov1@gmail.com,\\
a.yefimushkin@gmail.com

\end{document}